\documentclass[11pt]{article}

\usepackage[utf8]{inputenc}
\usepackage{amsmath}
\usepackage{amsfonts,amssymb,latexsym,multicol,color}
\usepackage[francais]{babel}
\usepackage[T1]{fontenc}

\newcounter{fig}
\newtheorem{theoa}{Th\'eor\`eme A}
\newtheorem{theob}{Th\'eor\`eme B}

\newtheorem{prop}{Proposition}

\newcommand{\cad}{\text{c'est-\`a-dire }}
\newcommand{\expli}[1]{\quad\text{\footnotesize (#1)}}

\makeatletter
\ifcase\@ptsize

\or

\or

\fi
\makeatother

\newcommand{\Implique}{\Longrightarrow}
\newcommand{\implique}{\Rightarrow}

\newcommand{\ioe}{\leqslant}
\newcommand{\soe}{\geqslant}

\newcommand{\vers}{\rightarrow}

\newcommand{\ssi}{\Leftrightarrow}

\newcommand{\demi}{{\frac{1}{2}}}

\newcommand{\Rcal}{{\mathcal R}}

\newcommand{\Nat}{{\mathbb N}}

\newcommand{\Real}{{\mathbb R}}

\newcommand{\fin}{\hfill$\Box$}
\newcommand{\dem}{\noindent {\bf D\'emonstration\ }}
\newcommand{\fine}{\tag*{\mbox{$\Box$}}}

\providecommand{\bysame}{\leavevmode ---\ }
\providecommand{\og}{``}
\providecommand{\fg}{''}
\providecommand{\smfandname}{et}
\providecommand{\smfedsname}{\'eds.}
\providecommand{\smfedname}{\'ed.}
\providecommand{\smfmastersthesisname}{M\'emoire}
\providecommand{\smfphdthesisname}{Th\`ese}

\newcommand{\virg}{\raisebox{.7mm}{,}}

\title{Sur la variation de certaines suites de parties fractionnaires}

\author{Michel Balazard, Leila Benferhat et Mihoub Bouderbala}

\date{}

\begin{document}
\maketitle

\begin{center}
  {\sc Abstract}
\end{center}
\begin{quote}
{\footnotesize Let $b>a>0$. We prove the following asymptotic formula
\begin{equation*}
\sum_{n\soe 0} \big\lvert\{x/(n+a)\}-\{x/(n+b)\}\big\rvert=\frac{2}{\pi}\zeta(3/2)\sqrt{cx}+O(c^{2/9}x^{4/9}),\\
\end{equation*}
with $c=b-a$, uniformly for $x \soe 40 c^{-5}(1+b)^{27/2}$.
 }
\end{quote}

\begin{center}
  {\sc Keywords}
\end{center}
\begin{quote}
{\footnotesize Fractional part, Elementary methods, van der Corput estimates }
\end{quote}



\section{Introduction}

Notons $\{t\}=t-\lfloor t\rfloor$ la partie fractionnaire du nombre réel $t$, où $\lfloor t \rfloor$ est la partie entière de~$t$. Pour $x>0$ et $b>a> 0$, les différences de parties fractionnaires 
\[
\{x/(n+a)\} - \{x/(n+b)\}\quad (n=0,1,2,\dots)
\]
sont les termes d'une série absolument convergente, puisque, pour $n>x-a$, cette différence vaut $cx/(n+a)(n+b)$, avec 
\[
c=b-a,
\]
notation que nous conserverons dans tout cet article. On peut donc considérer la norme au sens~$\ell^1$ de cette suite, \cad la quantité
\[
W(x;a,b)=\sum_{n\soe 0} \big\lvert\{x/(n+a)\}-\{x/(n+b)\}\big\rvert.
\]

La somme $W(x;1,2)$ joue un rôle auxiliaire dans l'article \cite{MR0016389} de Wintner. Il y démontra l'ordre de grandeur $W(x;1,2) \asymp \sqrt{x}$ (pour $x \soe 1$), et en déduisit l'optimalité de l'estimation~$O(\sqrt{x})$ pour le terme d'erreur de formules asymptotiques pour certaines moyennes arithmétiques. L'estimation de Wintner a été précisée par le premier auteur : on a
\begin{equation}\label{t0}
W(x;1,2)=\frac{2}{\pi}\zeta(3/2)\sqrt{x}+O(x^{2/5}) \quad (x>0),
\end{equation}
où $\zeta$ désigne la fonction zêta de Riemann (cf. \cite{zbMATH06754310}).

Le but du présent article est de généraliser \eqref{t0} à la somme $W(x;a,b)$. Afin d'énoncer notre premier résultat, il nous faut introduire les quantités
\begin{align}
R_0(x;a,b) &=\frac cx\sum_{1\ioe k\ioe K(x/c)} k^2\big(\{x/k-a\}-\{x/k-b\}\big)\label{181007a}\\
R_j(x;a,b) &=\sum_{K_{j-1}(x/c)<k\leqslant K_{j+1}(x/c)} (k^2c/x-j)\big(\{x/k-a\}-\{x/k-b\}\big) \quad (j \in \Nat^*),\label{181007b}
\end{align}
où, pour $t>0$ et $j \in \Nat$, nous notons 

$ \bullet \; K(t)$ le plus grand nombre entier $k$ tel que $k(k+1)\ioe t$ ;

$\bullet\; K_j(t)$ le plus grand nombre entier $k \soe j$ tel que~$(k-j)k\ioe jt$ (en particulier, $K_0(t) = 0$).

\smallskip

Enfin, nous posons, pour $J$ réel et positif,
\[
\Rcal(J,x;a,b) =\sum_{0\ioe j \ioe J} R_j(x;a,b).
\]

Observons que, pour $c$ entier, en particulier si $a=1$ et $b=2$, les quantités $R_j(x;a,b)$ sont nulles ; elles ne jouaient donc aucun rôle dans l'étude effectuée dans \cite{zbMATH06754310}. Notre premier résultat est une généralisation de \eqref{t0}.

\begin{theoa}
Pour $x \soe 40 c^{-3}(1+b)^4$, on a
\[
W(x;a,b)=\frac{2}{\pi}\zeta(3/2)\sqrt{cx}+\Rcal(J,x;a,b)+O\big((1+b)^{2/5}c^{1/5}x^{2/5}\big),
\]
où $J= c^{3/5}(1+b)^{-4/5}x^{1/5}$.
\end{theoa}


La somme $\Rcal(J,x;a,b)$ peut être estimée grâce aux résultats classiques de van der Corput, obtenus grâce à l'utilisation de sommes trigonométriques. Nous obtenons le résultat suivant.

\begin{theob}
Pour $x \soe 40 c^{-5}(1+b)^{27/2}$,
\[
W(x;a,b)=\frac{2}{\pi}\zeta(3/2)\sqrt{cx}+O(c^{2/9}x^{4/9}).
\]
\end{theob}

Observons que l'on en déduit l'estimation $W(x;a,b)\ll \sqrt{cx}$ sous la même hypothèse. 

\smallskip

La quantité $W(x;a,b)$ est reliée à la suivante, définie pour $x>0$ et $b>a>0$ par
\[
V(x;a,b)=\sum_{n\soe 0} \big(\{x/(n+a)\}-\{x/(n+b)\}\big).
\]
On a $V(x;1,2)=\{x\}$ et, plus généralement,
\[
V(x;a,b)=\{x/a\}+\cdots+\{x/(b-1)\}
\]
si $c=b-a$ est entier, mais l'estimation de $V(x;a,b)$ dans le cas général est un problème non trivial. 

Cette somme intervient dans l'étude de la question suivante. Soit $f$ une fonction arithmétique de période $q \in \Nat^*$, et de moyenne nulle. Sa fonction sommatoire $F$ est donc aussi périodique, de période $q$. Considérons le produit de convolution $g=f*{\mathbf 1}$. On a alors
\begin{equation*}
G(x)=\sum_{n\ioe x}g(n)=\sum_{n\soe 1} f(n)\lfloor x/n \rfloor = Cx -\Delta(x),
\end{equation*}
où $C=\sum_{n \soe 1} f(n)/n$ et
\begin{multline*}
\Delta(x)=\sum_{n\soe 1}f(n)\{x/n\}=\sum_{n\soe 1}F(n)\big(\{x/n\}-\{x/(n+1)\}\big)=\\
\sum_{k=1}^{q}F(k)\sum_{j\soe 0}\big(\{x/(jq+k)\}-\{x/(jq+k+1)\}\big)=\sum_{k=1}^{q}F(k)V\big(x/q;k/q,(k+1)/q\big).
\end{multline*}

La connaissance du comportement de $V(x;a,b)$ est donc susceptible d'apporter des informations sur celui du terme d'erreur $\Delta(x)$. La méthode de démonstration des théorèmes A et~B ci-dessus s'applique également à l'étude de la somme $V(x;a,b)$. Cela étant, la forme plus simple de cette quantité, relativement à $W(x;a,b)$, se prête \emph{a priori} à un traitement élémentaire classique via la méthode de l'hyperbole, suivi d'une application de la théorie de van der Corput, ou à une étude analytique à l'aide de la fonction $\zeta$ d'Hurwitz. Pour conserver au présent texte une unité méthodologique, nous n'y abordons donc pas l'étude de $V(x;a,b)$, nous contentant de signaler ici la majoration évidente $\lvert V(x;a,b)\rvert \ioe W(x;a,b)$. En particulier, en utilisant la majoration~$W(x;a,b)\ll \sqrt{cx}$, valable sous les conditions du Théorème B, on obtient l'estimation uniforme
\[
\Delta(x) \ll \frac{\sum_{k=1}^q\lvert F(k)\rvert}q\sqrt{x} \quad (1 \ioe q \ioe x^{1/6}/22).
\]
L'étude de $V(x;a,b)$ pourrait permettre de préciser ce résultat.

\smallskip

Le plan de cet article est le suivant. Au \S\ref{180827a}, nous décomposons $W(x;a,b)$ en somme de quantités $W_j(x;a,b)$, regroupant les entiers $n$ tels que
\[
\lfloor x/(n+a)\rfloor-\lfloor x/(n+b)\rfloor=j \quad ( j \in \Nat),
\]
et nous donnons des expressions de $W_j(x;a,b)$ (formules \eqref{t16} au \S\ref{180827c}, et \eqref{180515d} au \S\ref{181017a}). Au \S\ref{180827b}, nous donnons des estimations des quantités $W_j(x;a,b)$, faisant intervenir les sommes $R_j(x;a,b)$ définies par \eqref{181007a} et \eqref{181007b} ci-dessus. Nous en déduisons le Théorème A au \S\ref{180608c}. Enfin, au \S\ref{181020a}, nous exposons quelques éléments de la théorie de van der  Corput ; ils sont ensuite utilisés pour estimer les quantités $R_j(x;a,b)$. Cela nous permet d'obtenir le Théorème B.


\section{Décomposition de la somme $W(x;a,b)$}\label{180827a}

Comme $x$, $a$ et $b$ sont fixés dans ce paragraphe et les paragraphes \ref{180827b} et \ref{180608c}, nous allégeons la notation en écrivant simplement~$W$ au lieu de $W(x;a,b)$, et nous adopterons la même convention pour les quantités et ensembles, dépendant de $x$, $a$ et $b$, intervenant dans la démonstration. Les lettres $j,k,h,n$ désigneront toujours des variables entières positives ou nulles.

\subsection{Les sommes $W_j$}

En utilisant la notation d'Iverson ($[P]=1$ si la propriété $P$ est vraie, $[P]=0$ sinon), posons, pour $j\in \Nat$, 
$$
W_j=\sum_{n\soe 0} [\lfloor x/(n+a)\rfloor-\lfloor x/(n+b)\rfloor=j]\cdot\big\lvert\{x/(n+a)\}-\{x/(n+b)\}\big\rvert\, ,
$$
de sorte que
$$
W=\sum_{j\in\Nat}W_j.
$$

Nous allons évaluer $W_j$ en suivant la méthode adoptée dans \cite{zbMATH06754310}. Par souci de lisibilité, nous reproduisons, \emph{mutatis mutandis}, les détails des transformations opérées sur ces sommes.

\smallskip

Pour commencer, les relations
\begin{align*}
k&=\lfloor x/(n+a)\rfloor\\
h&=\lfloor x/(n+b)\rfloor
\end{align*}
entraînent $0\ioe h\ioe k \ioe x/a$ et équivalent à
\begin{align*}
k\ioe x/(n+a) <k+1\\
h\ioe x/(n+b)<h+1,
\end{align*}
autrement dit
\begin{equation}\label{t10}
\max\Big(\frac{x}{k+1}-a,\frac{x}{h+1}-b\Big)<n\ioe \min\Big(\frac{x}{k}-a,\frac{x}{h}-b\Big).
\end{equation}
(avec la convention $x/0=\infty$). 

\smallskip

Nous d\'esignerons par $I(h,k)$ l'intervalle de valeurs de $n$ d\'efini par l'encadrement \eqref{t10}, \cad
$$
I(h,k)=\Nat\,\cap\, ]x/(k+1)-a,x/k-a]\,\cap\,]x/(h+1)-b,x/h-b].
$$
La collection des $I(h,k)$ non vides constitue une partition de $\Nat$.

Notons que, pour $n\in I(h,k)$, on a
$$
\{x/(n+a)\}-\{x/(n+b)\}=cx/(n+a)(n+b)-k+h\virg
$$
où nous rappelons que $c$ désigne la différence $b-a$. En particulier,
\begin{equation*}
0 \ioe k-h < y+1,
\end{equation*}
où l'on a posé $y=cx/ab$.

La somme $W_j$ est donc nulle si $j\soe y+1$. Posons 
$$
W(h,k)=\sum_{n\in I(h,k)}\big\lvert cx/(n+a)(n+b)-k+h\big\rvert
$$
et, pour $j \in \Nat$ tel que $0\ioe j < y+1$, 
\begin{align}
E_j&=\{(h,k)\in\Nat^2,\, 0\ioe h\ioe k\ioe x/a, \, k-h=j\}\notag\\
&=\{(k-j,k),\, k\in \Nat, \, j \ioe k \ioe x/a\}\label{t13}
\end{align}
Nous aurons
$$
W_j=\sum_{(h,k)\in E_j}W(h,k).
$$

Nous allons établir une expression de la somme $W_j$, en commen\c{c}ant par le cas $j=0$.

\subsection{Expression de $W_0$}\label{180827c}

Nous avons 
$$
E_0=\{(k,k),\, k\in \Nat,\, 0\ioe k\ioe x/a\}.
$$
L'intervalle $I(k,k)$ de $\Nat$ est d\'efini par l'encadrement
\begin{equation}\label{t12}
\frac{x}{k+1}-a<n\ioe \frac{x}{k}-b.
\end{equation}

Il est vide si $k(k+1) > x/c$. Pour $t>0$, d\'esignons par $K(t)$ le plus grand nombre entier $k$ tel que $k(k+1)\ioe t$, et observons simplement, pour l'instant, que $K(t) \ioe \sqrt{t}$. 

Posons également
\begin{equation}\label{t3}
F(t)= \sum_{n > t}1/(n+a)(n+b) \quad (t\soe 0)
\end{equation}
et, par convention, 
\begin{align*}
F(\infty)&=0\\
F(t)=F(0^-)&=\sum_{n \soe 0}1/(n+a)(n+b) \quad(t<0).
\end{align*}

Pour $x\soe a^2/c$, en notant $K=K(x/c)$, on aura $K \ioe \sqrt{x/c}\ioe x/a$, donc
\begin{align}
W_0 &=\sum_{0\ioe k\ioe K}W(k,k)\notag\\
&=\sum_{0\ioe k\ioe K}\quad\sum_{x/(k+1)-a<n\ioe x/k-b}cx/(n+a)(n+b)\notag\\
&=cx\sum_{0\ioe k\ioe K} \Big(F\big( x/(k+1)-a\big) -F(x/k-b )\Big)\notag\\
&=cxF\big(x/(K+1)-a\big) + cx\sum_{1\ioe k\ioe K} \big(F(x/k-a) -F(x/k-b)\big).\label{t16}
\end{align}

\subsection{D\'ecomposition de l'ensemble $E_j$}

Nous supposons maintenant $j \soe 1$. Toujours en adaptant la démarche suivie dans \cite{zbMATH06754310}, nous allons d\'ecomposer l'ensemble $E_j$ d\'efini par \eqref{t13}
en une partition de trois sous-ensembles sur lesquels l'encadrement \eqref{t10} s'exprimera sans recours aux fonctions $\max$ et $\min$.

\smallskip

Si $k>h\soe 0$ et $x>0$, l'in\'egalit\'e
$$
\frac xk-a \ioe \frac xh -b
$$
\'equivaut \`a
$$
\frac{hk}{k-h} \ioe x/c.
$$
En particulier, on  a les implications
$$
\frac{x}{k+1}-a\ioe\frac{x}{h+1}-b \Implique \frac xk-a \ioe \frac xh -b
$$
et
$$
\frac xk-a > \frac xh -b \Implique \frac{x}{k+1}-a > \frac{x}{h+1}-b.
$$

Nous consid\'erons donc les trois parties suivantes de $E_j$  (la d\'efinition de chaque $E_{j,i}$ est suivie par la forme que prend l'encadrement \eqref{t10} lorsque $(k-j,k)\in E_{j,i}$) :
\begin{align}
E_{j,1}&=\{(k-j,k) : j\ioe k\ioe x/a, \, (k-j+1)(k+1) \ioe jx/c\}\notag\\
&\quad \quad\quad \frac{x}{k-j+1}-b<n\ioe \frac{x}{k}-a\label{t8}\\ \notag\\
E_{j,2}&=\{(k-j,k) : j\ioe k\ioe x/a, \, (k-j)k > jx/c\}\notag\\
&\quad\quad\quad\frac{x}{k+1}-a<n\ioe \frac{x}{k-j}-b\label{t9}\\ \notag\\
E_{j,3}&=\{(k-j,k) : j\ioe k\ioe x/a, \, (k-j)k\ioe jx/c <(k-j+1)(k+1)\}\notag\\
&\quad\quad\quad\frac{x}{k+1}-a<n\ioe \frac{x}{k}-a\label{t5}
\end{align}

\smallskip

Celles des trois parties $E_{j,i}$ ($1\ioe i\ioe 3$) qui sont non vides forment une partition de $E_j$. Par cons\'equent, on a
$$
W_j=W_{j,1}+W_{j,2}+W_{j,3},
$$
o\`u
$$
W_{j,i}=\sum_{(h,k)\in E_{j,i}} W(h,k) \quad (1\ioe i \ioe 3).
$$

Avant d'\'evaluer successivement les trois quantit\'es $W_{j,i}$, nous allons définir, aux sous-paragraphes suivants, deux fonctions auxiliaires, $K_j$ et $N_j$.

\subsection{La fonction $K_j(t)$}\label{180528a}

Pour $j \in \Nat$ et $t > 0$, nous définissons $K_j(t)$ comme le plus grand nombre entier $k \soe j$ tel que~$(k-j)k\ioe jt$. En particulier, $K_0(t) = 0$. 


\smallskip

Au moyen de la fonction $K_j(t)$, on peut, pour $j \soe 1$, r\'ecrire les conditions, quadratiques relativement \`a~$k$, intervenant dans les d\'efinitions des ensembles $E_{j,i}$, sous les formes suivantes, respectivement :
\begin{align}
k&\ioe K_j(x/c)-1 & (i=1)\label{u16}\\
k&>K_j(x/c) & (i=2)\label{t17}\\
k&=K_j(x/c) & (i=3)\label{t18}
\end{align}

De plus, la condition $k\ioe x/a$, qui figure également dans la définition de ces ensembles, est superflue pour $i=1,3$, si $j \ioe y=cx/ab$. En effet la relation $j\ioe y$ peut s'écrire sous la forme
\[
j x/c \ioe \Big(\frac x{a}-j\Big)\frac x{a}\cdotp
\]
Les inégalités 
\[
\big(K_j(x/c)-j\big)K_j(x/c)\ioe jx/c \quad ; \quad j\ioe y\ioe x/a, 
\]
et le fait que $t \mapsto (t-j)t$ est strictement croissante pour $t \soe j$ entraînent alors l'inégalité 
\begin{equation}\label{180508a}
K_j(x/c)\ioe x/a.
\end{equation}

\subsection{La fonction $N_j(x;a,b)$}

Pour exprimer la quantit\'e $\lvert j-cx/(n+a)(n+b)\rvert$ sans valeur absolue, nous sommes conduits \`a d\'efinir, pour $j\in \Nat^*$, $N_j=N_j(x;a,b)$ comme le plus grand nombre entier $n$ tel que
$$
(n+a)(n+b)\ioe cx/j.
$$
On a donc $N_j\soe 0$ si $j\ioe y=cx/ab$.

\smallskip

\'Etablissons maintenant une relation entre les quantités $K_j=K_j(x/c)$ et $N_j=N_j(x;a,b)$.
 
\begin{prop}\label{t14}
Si $j \in \Nat^*$ et $j \ioe y$, on a
$$
\Big\lfloor \frac x{K_j+1}-a\Big\rfloor \ioe N_j \ioe \Big\lfloor \frac x{K_j}-a\Big\rfloor.
$$
\end{prop}
\dem

L'encadrement définissant $K_j$,
$$
(K_j-j)K_j \ioe jx/c < (K_j+1-j)(K_j+1)\, ,
$$
équivaut à
$$
\frac{x}{K_j+1}\Big( \frac{x}{K_j+1} +c\Big) < \frac{cx}j \ioe \frac{x}{K_j}\Big( \frac{x}{K_j} +c\Big),
$$
\cad à
\[
\Big(\frac{x}{K_{j+1}}-a+a\Big)\cdot\Big( \frac{x}{K_{j+1}}-a +b\Big)< \frac{cx}j \ioe\Big(\frac{x}{K_j}-a+a\Big)\cdot\Big( \frac{x}{K_j}-a +b\Big).
\]

On a $N_j \soe 0$, et $t\mapsto (t+a)(t+b)$ est strictement croissante pour $t\soe -a$. Le dernier encadrement entraîne donc celui de l'énoncé.\fin

\subsection{Calcul de $W_{j,1}$}

Supposons $1 \ioe j \ioe y$. En notant simplement $K_j$ pour $K_j(x/c)$, nous aurons, d'après \eqref{t8},~\eqref{u16} et \eqref{180508a},
\begin{align*}
W_{j,1} &=\sum_{(h,k)\in E_{j,1}} W(h,k)\\
&=\sum_{j \ioe k\ioe K_j-1}\;\sum_{\frac{x}{k-j+1}-b<n\ioe \frac{x}{k}-a}\big\lvert cx/(n+a)(n+b)-j\big\rvert.
\end{align*}

La somme int\'{e}rieure est non vide seulement si
\begin{equation*}
\frac{x}{k-j+1}-b<\frac{x}{k}-a,
\end{equation*}
autrement dit seulement si $k>K_{j-1}$(rappelons que $K_{0}=0).$ Les nombres
entiers $n$ intervenant dans cette somme int\'{e}rieure sont strictement sup\'{e}rieurs \`{a} 
\begin{equation*}
\frac{x}{K_{j}-j}-b\geqslant \frac{x}{K_{j}}-a\geqslant N_{j},
\end{equation*}
d'apr\`{e}s la proposition \ref{t14}.

Dans le calcul qui suit, ainsi qu'au paragraphe suivant, nous emploierons la
fonction $F$ définie par \eqref{t3}, et l'identit\'{e} 
\guillemotleft $r-$t\'{e}lescopique\guillemotright,
\begin{equation*}
\sum_{\alpha <k\leqslant \beta }\left( u_{k}-u_{k-r}\right) =\sum_{\beta
-r<k\leqslant \beta }u_{k}-\sum_{\alpha -r<k\leqslant \alpha }u_{k}.
\end{equation*}

On a donc
\begin{align}
W_{j,1} &=\sum_{K_{j-1}<k\leqslant K_{j}-1}\sum_{\frac{x}{k-j+1}
-b<n\leqslant \frac{x}{k}-a}\big( j-cx/(n+a)( n+b) \big) \notag\\
&=j\sum_{K_{j-1}<k\leqslant K_{j}-1}\left( \left\lfloor \frac{x}{k}
-a\right\rfloor -\left\lfloor \frac{x}{k-j+1}-b\right\rfloor \right)\notag \\
&\quad-cx\sum_{K_{j-1}<k\leqslant K_{j}-1}\left( F\left( \frac{x}{k-j+1}-b\right)
-F\left( \frac{x}{k}-a\right) \right) .\label{180508b}
\end{align}

L'avant-dernière somme vaut
\begin{multline*}
\sum_{K_{j-1}<k\leqslant K_{j}-1}\left( \left\lfloor \frac{x}{k}
-a\right\rfloor -\left\lfloor \frac{x}{k-j+1}-b\right\rfloor \right) =\\
\sum_{K_{j-1}<k\leqslant K_{j}-1}\left( \left\lfloor \frac{x}{k}-a\right\rfloor -  \left\lfloor \frac{x}{k} -b\right\rfloor   \right)+\sum_{K_{j-1}<k\leqslant K_{j}-1}\left( \left\lfloor \frac{x}{k}-b\right\rfloor -\left\lfloor \frac{x}{k-j+1}-b\right\rfloor \right),
\end{multline*}
où, par l'identité {\og$(j-1)$-télescopique\fg},
\begin{multline*}
\sum_{K_{j-1}<k\leqslant K_{j}-1}\left( \left\lfloor \frac{x}{k}-b\right\rfloor -\left\lfloor \frac{x}{k-j+1}-b\right\rfloor \right)=\\
\sum_{K_{j}-j<k\leqslant K_{j}-1} \left\lfloor \frac{x}{k}-b\right\rfloor-\sum_{K_{j-1}-(j-1)<k\leqslant K_{j-1}} \left\lfloor \frac{x}{k}-b\right\rfloor.
\end{multline*}

Une manipulation similaire s'applique à la dernière somme de \eqref{180508b}, et on obtient finalement
\begin{multline}\label{180515a}
W_{j,1} =j\sum_{K_{j-1}<k\leqslant K_{j}-1}\left( \left\lfloor \frac{x}{k}-a\right\rfloor -  \left\lfloor \frac{x}{k} -b\right\rfloor   \right)\\
+j\sum_{K_{j}-j<k\leqslant K_{j}-1} \left\lfloor \frac{x}{k}-b\right\rfloor-j\sum_{K_{j-1}-(j-1)<k\leqslant K_{j-1}} \left\lfloor \frac{x}{k}-b\right\rfloor\\
+cx\sum_{K_{j-1}<k\leqslant K_{j}-1}\left( F\left( \frac{x}{k}-a\right) -  F\left( \frac{x}{k} -b\right)   \right)\\
+cx\sum_{K_{j}-j<k\leqslant K_{j}-1} F \left( \frac{x}{k}-b\right)-cx\sum_{K_{j-1}-(j-1)<k\leqslant K_{j-1}} F\left( \frac{x}{k}-b\right).
\end{multline}

\subsection{Calcul de $W_{j,2}$}

On a ici, toujours pour $1 \ioe j \ioe y$, et d'après \eqref{t9},~\eqref{t17} et \eqref{180508a},
\begin{align*}
W_{j,2} &=\sum_{(h,k)\in E_{j,2}} W(h,k)\\
&=\sum_{K_j < k\ioe x/a}\;\sum_{\frac{x}{k+1}-a<n\ioe \frac{x}{k-j}-b}\big\lvert cx/(n+a)(n+b)-j\big\rvert.
\end{align*}

La somme int\'{e}rieure est non vide seulement si
\begin{equation*}
\frac{x}{k+1}-a\ioe \frac{x}{k-j}-b,
\end{equation*}
autrement dit seulement si $k\ioe K_{j+1}-1$. Nous supposerons donc maintenant que $j\ioe y-1$, de sorte que $K_{j+1} \ioe x/a$, d'après \eqref{180508a}. On obtient alors
\[
W_{j,2}=\sum_{K_j < k\ioe K_{j+1}-1}\;\sum_{\frac{x}{k+1}-a<n\ioe \frac{x}{k-j}-b}\big\lvert cx/(n+a)(n+b)-j\big\rvert.
\]

Les nombres entiers $n$ intervenant dans la somme int\'{e}rieure sont inf\'{e}rieurs ou égaux \`{a} 
\begin{equation*}
\frac{x}{K_{j}+1-j}-b < \frac{x}{K_{j}+1}-a,
\end{equation*}
par définition de $K_j$. La proposition \ref{t14} prouve alors que ces nombres entiers sont $\ioe N_j$.

On a donc
\begin{align}
W_{j,2} &=\sum_{K_j < k\ioe K_{j+1}-1}\;\sum_{\frac{x}{k+1}-a<n\ioe \frac{x}{k-j}-b}\big ( cx/(n+a)(n+b)-j\big)\notag\\
&=j\sum_{K_j < k\ioe K_{j+1}-1}\left(\left\lfloor \frac{x}{k+1}-a\right\rfloor - \left\lfloor \frac{x}{k-j}
-b\right\rfloor \right)\notag \\
&\quad+cx\sum_{K_j < k\ioe K_{j+1}-1}\left( F\left( \frac{x}{k+1}-a\right)
-F\left( \frac{x}{k-j}-b\right) \right) .\label{180509b}
\end{align}

L'avant-dernière somme vaut
\begin{multline*}
\sum_{K_j +1< k\ioe K_{j+1}}\left( \left\lfloor \frac{x}{k}-a\right\rfloor-\left\lfloor \frac{x}{k-j-1}
-b\right\rfloor  \right) =\\
\sum_{K_{j}+1<k\leqslant K_{j+1}}\left( \left\lfloor \frac{x}{k}-a\right\rfloor -  \left\lfloor \frac{x}{k} -b\right\rfloor   \right)+\sum_{K_{j}+1<k\leqslant K_{j+1}}\left( \left\lfloor \frac{x}{k}-b\right\rfloor -\left\lfloor \frac{x}{k-j-1}-b\right\rfloor \right),
\end{multline*}
où, par l'identité {\og$(j+1)$-télescopique\fg},
\begin{multline*}
\sum_{K_{j}+1<k\leqslant K_{j+1}}\left( \left\lfloor \frac{x}{k}-b\right\rfloor -\left\lfloor \frac{x}{k-j-1}-b\right\rfloor \right)=\\
\sum_{K_{j+1}-j-1<k\leqslant K_{j+1}} \left\lfloor \frac{x}{k}-b\right\rfloor-\sum_{K_{j}-j<k\leqslant K_{j}+1} \left\lfloor \frac{x}{k}-b\right\rfloor.
\end{multline*}

Une manipulation similaire s'applique à la somme \eqref{180509b}, et on obtient finalement
\begin{multline}\label{180515b}
W_{j,2} =j\sum_{K_{j}+1<k\leqslant K_{j+1}}\left( \left\lfloor \frac{x}{k}-a\right\rfloor -  \left\lfloor \frac{x}{k} -b\right\rfloor   \right)\\
+j\sum_{K_{j+1}-j-1<k\leqslant K_{j+1}} \left\lfloor \frac{x}{k}-b\right\rfloor-j\sum_{K_{j}-j<k\leqslant K_{j}+1} \left\lfloor \frac{x}{k}-b\right\rfloor\\
+cx\sum_{K_{j}+1<k\leqslant K_{j+1}}\left( F\left( \frac{x}{k}-a\right) -  F\left( \frac{x}{k} -b\right)   \right)\\
+cx\sum_{K_{j+1}-j-1<k\leqslant K_{j+1}}F \left( \frac{x}{k}-b\right)-cx\sum_{K_{j}-j<k\leqslant K_{j}+1} F\left( \frac{x}{k}-b\right).
\end{multline}

\subsection{Calcul de $W_{j,3}$}

Pour $1\leqslant j\leqslant y$, on a, d'apr\`{e}s \eqref{t5}, \eqref{t18}, et \eqref{180508a},

\begin{align*}
W_{j,3} &=\sum_{\frac{x}{K_{j}+1}-a<n\leqslant \frac{x}{K_{j}}
-a}\left\vert j-cx/(n+a)\left( n+b\right) \right\vert \\
&=\sum_{\frac{x}{K_{j}+1}-a<n\leqslant N_{j}}\big( cx/(n+a)\left(
n+b\right) -j\big) +\sum_{N_{j}<n\leqslant \frac{x}{K_{j}}-a}\big(
j-cx/(n+a)\left( n+b\right) \big) ,
\end{align*}
d'apr\`{e}s la proposition \ref{t14}. Observons que l'avant-derni\`{e}re somme est vide si $N_{j}=\left\lfloor 
\frac{x}{K_{j}+1}-a\right\rfloor .$

On a 
\begin{align*}
\sum_{\frac{x}{K_{j}+1}-a<n\leqslant N_{j}}\big( cx/(n+a)\left( n+b\right)-j\big) &=cx\left( F\left( \frac{x}{K_{j}+1}-a\right) -F\left(N_{j}\right) \right)\\
& \qquad \qquad\qquad\qquad -j\left( N_{j}-\left\lfloor \frac{x}{K_{j}+1}
-a\right\rfloor \right) \\
\sum_{N_{j}<n\leqslant \frac{x}{K_{j}}-a}\big( j-cx/(n+a)\left( n+b\right)
\big) &=j\left( \left\lfloor \frac{x}{K_{j}}-a\right\rfloor -N_{j}\right)
-cx\left( F\left( N_{j}\right) -F\left(\frac{x}{K_{j}}-a\right)\right) ,
\end{align*}
donc
\begin{multline}\label{180515c}
W_{j,3} =cx\left(F\left( \frac{x}{K_{j}+1}-a\right) +F\left(\frac{x}{K_{j}}
-a\right)-2F\left( N_{j}\right) \right)  \\
+j\left(\left\lfloor \frac{x}{K_{j}}-a\right\rfloor + \left\lfloor \frac{x}{K_{j}+1}-a\right\rfloor -2N_{j}
\right) .  
\end{multline}

\subsection{Expression de $W_j$ pour $j \soe 1$}\label{181017a}

En regroupant les identités \eqref{180515a}, \eqref{180515b} et \eqref{180515c}, on obtient, pour $1 \ioe j \ioe y-1$,
\begin{multline}\label{180515d}
W_j =j\sum_{K_{j-1}<k\leqslant K_{j+1}}\left( \left\lfloor \frac{x}{k}-a\right\rfloor -  \left\lfloor \frac{x}{k} -b\right\rfloor   \right)\\
+j\sum_{K_{j+1}-(j+1)<k\leqslant K_{j+1}} \left\lfloor \frac{x}{k}-b\right\rfloor-j\sum_{K_{j-1}-(j-1)<k\leqslant K_{j-1}} \left\lfloor \frac{x}{k}-b\right\rfloor -2jN_j\\
+cx\sum_{K_{j-1}<k\leqslant K_{j+1}}\left( F\left( \frac{x}{k}-a\right) -  F\left( \frac{x}{k} -b\right)   \right)\\
+cx\sum_{K_{j+1}-(j+1)<k\leqslant K_{j+1}} F \left( \frac{x}{k}-b\right)-cx\sum_{K_{j-1}-(j-1)<k\leqslant K_{j-1}} F\left( \frac{x}{k}-b\right)-2cxF(N_j).
\end{multline}

\section{Estimation des sommes $W_j$}\label{180827b}

Nous allons maintenant utiliser les identités obtenues au paragraphe précédent pour estimer les contributions à $W$ des sommes $W_j$. Nous commençons par le cas $j=0$.

\subsection{La fonction $K$}

Au \S\ref{180827c}, pour $t>0$, nous avons introduit la notation $K(t)$ pour désigner le plus grand nombre entier $k$ tel que~$k(k+~1)\ioe ~t$, \cad
$$
K(t)=\lfloor\sqrt{t+1/4}-1/2\rfloor,
$$
et noté simplement $K$ la valeur $K(x/c)$. Nous utiliserons l'encadrement
\[
\sqrt{t}-3/2\ioe \sqrt{t+1/4}-3/2\ioe K(t)\ioe \sqrt{t}.
\]
\begin{prop}\label{180518a}
Pour $x\soe 9b^2/c$, on a 
\[
\frac x{K+1}-b \soe \frac{\sqrt{cx}}6\cdotp
\]
\end{prop}
\dem

On a
\begin{equation*}
\frac x{K+1}-b \soe cK-b \soe c(\sqrt{x/c}-3/2)-b \soe \sqrt{cx}-5b/2 \soe \frac{\sqrt{cx}}6\virg
\end{equation*}
si $x\soe 9b^2/c$.\fin

\subsection{Estimation de la fonction $F$}

Rappelons la définition \eqref{t3} :
\begin{equation*}
F(t)=\sum_{n >t}1/(n+a)(n+b).
\end{equation*}

\begin{prop}
Pour $b>a>0$ et $t>0$, on a
\begin{align}
F(t)&=\frac 1t +  O\left((1+b)/t^2\right).\label{180519a}\\
F(t)&=\frac 1t +  \frac{\{t\}-(a+b+1)/2}{t^2}+O\left((1+b)^2/t^3\right)\label{180518b}
\end{align}
\end{prop}
\dem

Nous démontrons \eqref{180518b} ; la démonstration de \eqref{180519a} est similaire, et plus simple.

Pour $n \in \Nat^*$, on a
\[
\frac 1{(n+a)(n+b)}-\frac{1}{n^2}+\frac{a+b}{n^3}=\frac{(a^2+ab+b^2)n+ab(a+b)}{n^3(n+a)(n+b)} \ioe \frac{3b^2}{n^4}\cdotp
\]

Pour $t>0$, on a
\begin{align*}
\sum_{n>t} \frac 1{n^2} &= \frac 1t + \frac{\{t\}-1/2}{t^2} +O(t^{-3})\\
\sum_{n>t} \frac 1{n^3} &= \frac 1{2t^2}  +O(t^{-3})\\
\sum_{n>t} \frac 1{n^4} &= O(t^{-3}).
\end{align*}

Par conséquent,
\begin{align*}
F(t) &=\sum_{n>t}\left(\frac{1}{n^2}-\frac{a+b}{n^3} + O(b^2/n^4)\right)\\
&=\frac 1t +  \frac{\{t\}-(a+b+1)/2}{t^2}+O\left((1+b)^2/t^3\right).\fine
\end{align*}

\subsection{Estimation de $W_0$}

Nous commençons par estimer le premier terme de l'expression \eqref{t16} de $W_0$.
\begin{prop}\label{180521b}
Pour $x\soe 9b^2/c$, on a
\[
F\big(x/(K+1)-a\big) =\frac 1{\sqrt{cx}}+O\big((1+b)/cx\big).
\]
\end{prop}
\dem

On a $K+1=\sqrt{x/c}+\vartheta_0$, où $\lvert \vartheta_0\rvert \ioe 1$ et $\sqrt{x/c} \soe 3b/c\soe 3$. Par conséquent,
\[
\frac x{K+1}-a=\frac x{\sqrt{x/c}+\vartheta_0} -a=\sqrt{cx}+O(b).
\]

On a donc $x/(K+1)-a=\sqrt{cx}(1+\vartheta_1)$, avec $\vartheta_1\soe -5/6$ (d'après la proposition \ref{180518a}), et~$\vartheta_1=O(b/\sqrt{cx})$. L'estimation \eqref{180519a} nous donne alors
\begin{equation*}
F\big(x/(K+1)-a\big) =\frac{1+O(b/\sqrt{cx})}{\sqrt{cx}} +O\big((1+b)/cx\big)=\frac 1{\sqrt{cx}}+O\big((1+b)/cx\big).\fine
\end{equation*}

\smallskip

Passons à la somme apparaissant dans \eqref{t16}.
\begin{prop}\label{180521a}
Pour $x\soe 9b^2/c$, on a
\begin{multline*}
\sum_{1\ioe k\ioe K} \big(F(x/k-a) -F(x/k-b)\big) =\\
-\frac 1{3\sqrt{cx}}+x^{-2}\sum_{1\ioe k\ioe K} k^2\big(\{x/k-a\}-\{x/k-b\}\big)+O\big((1+b)^2/c^2x\big).
\end{multline*}
\end{prop}
\dem

Nous allons utiliser l'estimation \eqref{180518b}. Par la proposition \ref{180518a}, les quantités $x/k-a$ et $x/k-b$, figurant dans la somme à évaluer, sont toutes $\soe \sqrt{cx}/6$. Par conséquent, la contribution à cette somme du terme d'erreur de \eqref{180518b} est
\[
\ll (1+b)^2K/(cx)^{3/2} \ioe (1+b)^2/c^2x.
\]

Pour $1\ioe k \ioe K$, on a
\begin{align*}
\frac 1{x/k-a}-\frac 1{x/k-b}&=-\frac{ck^2}{x^2}\left(1-\frac{ak}x\right)^{-1}\left(1-\frac{bk}x\right)^{-1}\\
&=-ck^2/x^2+O(cbk^3/x^3),
\end{align*}
car $bk/x \ioe 1/3$ si $1\ioe k \ioe K$.

Par conséquent, la contribution à la somme étudiée du terme $1/t$ de \eqref{180518b} vaut
\begin{align*}
\sum_{1\ioe k\ioe K}\left (\frac 1{x/k-a}-\frac 1{x/k-b}\right)&=-\frac{c}{x^2}\big(K^3/3+O(K^2)\big)+O(cbK^4/x^3)\\
&=-\frac 1{3\sqrt{cx}}+O(b/cx).
\end{align*}

Enfin,
\begin{multline*}
\frac{\{x/k-a\}-(a+b+1)/2}{(x/k-a)^2}-\frac{\{x/k-b\}-(a+b+1)/2}{(x/k-b)^2}=\\
\big(\{x/k-a\}-(a+b+1)/2\big)\big(k^2/x^2+O(bk^3/x^3)\big)\\
\qquad\qquad-\big(\{x/k-b\}-(a+b+1)/2\big)\big(k^2/x^2+O(bk^3/x^3)\big)\\
=k^2\big(\{x/k-a\}-\{x/k-b\}\big)/x^2+O\big((1+b)^2k^3/x^3\big).
\end{multline*}

La contribution du dernier terme d'erreur à la somme figurant dans \eqref{t16} est 
\[
\ll (1+b)^2K^4/x^3\ll (1+b)^2/c^2x.
\]

Le résultat découle de ces estimations.\fin

\smallskip

Les propositions \ref{180521b} et \ref{180521a} et la formule \eqref{t16} fournissent l'expression suivante de $W_0$.

\begin{prop}\label{180528j}
Pour $x\soe 9b^2/c$, on a
\[
W_0 =\frac 2{3}\sqrt{cx}+R_0+O\big((1+b)^2/c\big).
\]
où
\begin{equation}\label{180715a}
R_0=\frac cx\sum_{1\ioe k\ioe K} k^2\big(\{x/k-a\}-\{x/k-b\}\big)\quad (x>0).
\end{equation}
\end{prop}

\subsection{Grandes valeurs de $j$}

Avant d'examiner en détail chaque quantité $W_j$, notons l'estimation suivante, qui nous permettra de limiter les valeurs de $j$ à considérer.
\begin{prop}\label{180527d}
Soit $J$ un nombre r\'eel sup\'erieur \`a $1$. Pour $x>0$, et $b>a>0$, on a
\[
\sum_{j>J} W_j <\sqrt{cx/(J-1)}+1.
\]
\end{prop}
\dem

Si $j=\lfloor x/(n+a)\rfloor-\lfloor x/(n+b)\rfloor >J$, alors
\begin{equation*}
\frac{cx}{n^2} >\frac{cx}{(n+a)(n+b)}=j +\{ x/(n+a)\}-\{ x/(n+b)\}> J-1,
\end{equation*}
donc $n<\sqrt{cx/(J-1)}$. On en d\'eduit que
\begin{equation*}
\sum_{j>J}  W_j \ioe\sum_{0\ioe n<\sqrt{cx/(J-1)}}1<\sqrt{cx/(J-1)}+1.
\end{equation*}

\subsection{Résultats auxiliaires sur les quantités $K_j$}

Nous allons utiliser les résultats, démontrés dans \cite{zbMATH06754310}, concernant les fonctions $K_j(t)$, dont la définition a été rappelée au \S\ref{180528a}, ici évaluées en $t=x/c$. 

Pour $0\ioe j\ioe x/c$ et $x\soe c$, les propositions 2 et 3, p. 12 de \cite{zbMATH06754310}, affirment que 
\begin{align}
K_{j+1}-K_{j}&=\big(\sqrt{j+1}-\sqrt{j}\, \big)\sqrt{x/c}+O(1)\label{180522a}\\
\sum_{K_{j}<k\ioe K_{j+1}}k^2&=\frac{(j+1)\sqrt{j+1}-j\sqrt{j}}{3}\,(x/c)^{3/2}+O\big((j+1)x/c\big)\label{180522b}.
\end{align}

Notons que \eqref{180522a} entraîne
\begin{equation}\label{180528c}
K_{j+1}-K_{j} \ll (x/cj)^{1/2} \quad (0 <j \ioe x/c).
\end{equation} 

L'encadrement (17), p. 11 de \cite{zbMATH06754310} se récrit
\begin{equation}\label{180522c}
\sqrt{jx/c}+j/2-1< K_j\ioe \sqrt{jx/c}+j.
\end{equation} 
En particulier, pour $0 <j \ioe x/c$, on a
\begin{equation}\label{180522z}
 K_j\soe \demi\sqrt{jx/c}\,,
\end{equation} 
et
\begin{equation}\label{180711b}
K_{j+1} \ioe \sqrt{(j+1)x/c}+j+1 \ioe 2\sqrt{jx/c}+2j\ioe 4\sqrt{jx/c}.
\end{equation}

\smallskip

Les quatre propositions suivantes sont des lemmes utilisés lors des calculs des paragraphes suivants.

\begin{prop}\label{180711a}
Pour $0 <j \ioe x/c$ et $K_{j-1}< k \ioe K_{j+1}$, on a
\begin{equation*}
-1 < ck^2/x-j \ioe 1+8j^{3/2}(c/x)^{1/2}.
\end{equation*}
\end{prop}
\dem

En effet, par définition de $K_{j \pm 1}$, on a, pour $K_{j-1}< k \ioe K_{j+1}$,
\[
(j-1)x/c+(j-1)k < k^2 \ioe (j+1)x/c + (j+1)k,
\]
donc
\[
-1+(j-1)kc/x <ck^2/x-j \ioe 1+(j+1)kc/x.
\]

L'encadrement annoncé résulte alors de la majoration $K_{j+1} \ioe 4\sqrt{jx/c}$.\fin

\begin{prop}\label{180528d}
Pour $0 <j \ioe x/c$ 
\[
\sum_{K_{j-1}< k \ioe K_{j+1}}\lvert k^2c/x-j\rvert \ll(x/cj)^{1/2} +j.
\]
\end{prop}
\dem

En utilisant \eqref{180528c} et la proposition \ref{180711a}, on obtient
\begin{align*}
\sum_{K_{j-1}< k \ioe K_{j+1}}\lvert k^2c/x-j\rvert &\ll(x/cj)^{1/2} \big(1+j^{3/2}(c/x)^{1/2}\big)\\
&=(x/cj)^{1/2} +j.\fine
\end{align*}

\begin{prop}\label{180526i}
Pour $0<j\ioe cx/9(1+b)^2$ et $k\ioe K_j$, on a
\begin{align}
\lfloor x/k -b\rfloor &\soe \frac{5}{12}\sqrt{cx/j}\label{180526c}\\
x/k -b &\soe x/2k.\label{180526j}
\end{align}
\end{prop}
\dem

En utilisant \eqref{180522c}, on a, d'une part,
\begin{equation*}
\lfloor x/k -b\rfloor \soe x/K_j-b-1\soe \frac{x}{\sqrt{jx/c}+j}-b-1 \soe \frac 34\sqrt{cx/j}-b-1 \soe \frac{5}{12}\sqrt{cx/j}.
\end{equation*}

D'autre part,
\begin{equation*}
kb \ioe K_jb\ioe b\sqrt{jx/c}+jb \ioe \frac{bx}{3(1+b)}+\frac{bcx}{9(1+b)^2}\ioe \frac x2.\fine
\end{equation*}

\begin{prop}\label{180522d}
Pour $0<j\ioe cx/9(1+b)^2$ et $K_j-j<k\ioe K_j$, on a
$$
j\lfloor x/k -b\rfloor +cx/\lfloor x/k-b \rfloor=2\sqrt{jcx}+O\big((1+b)^2j^{3/2}(cx)^{-1/2}\big).
$$
\end{prop}
\dem

Il s'agit d'une adaptation de la proposition 5, p. 17 de \cite{zbMATH06754310}.

Posons $q=\lfloor x/k-b \rfloor$. On a
$$
jq+cx/q-2\sqrt{jcx}=p^2,
$$
o\`u
\begin{equation*}
p= \sqrt{jq}-\sqrt{cx/q}=\frac{jq^2-cx}{q\big(\sqrt{jq}+\sqrt{cx/q}\big)}.
\end{equation*}

D'après \eqref{180522c}, on a
\[
\sqrt{jx/c}+j \soe K_j \soe k \soe K_j-j+1\soe \sqrt{jx/c}-j/2\soe \frac 56\sqrt{jx/c},
\]
donc, compte également tenu de la proposition \ref{180526i},
\begin{equation}\label{t45}
q=\lfloor x/k-b \rfloor=x/k+O(1+b)=\frac{x}{\sqrt{jx/c}+O(j)}+O(1+b)=\sqrt{cx/j}+O(1+b).
\end{equation}

Par cons\'equent,
\begin{align*}
jq^2-cx&=O\big((1+b)\sqrt{cjx}\,\big)\\
p&=O\big((1+b)j^{3/4}(cx)^{-1/4}\big)
\end{align*}
et 
\begin{equation*}
jq+cx/q-2\sqrt{jcx}=O\big((1+b)^2j^{3/2}(cx)^{-1/2}\big).\fine
\end{equation*}

\subsection{Estimation de $W_j$ pour $j \soe 1$}

Commençons en observant que la condition $j\ioe cx/9(1+b)^2 -1$, qui va figurer dans les trois premiers énoncés de ce sous-paragraphe, entraîne l'inégalité $j \ioe y-1$, condition d'application de l'identité \eqref{180515d}.

Nous estimons en premier lieu la contribution à $W_j$ des sommes de \eqref{180515d} où $a$ ne figure pas.
\begin{prop}\label{180527c}
Pour $0<j\ioe cx/9(1+b)^2 -1$, on a
\begin{multline}\label{180526a}
j\sum_{K_{j+1}-(j+1)<k\leqslant K_{j+1}} \left\lfloor \frac{x}{k}-b\right\rfloor-j\sum_{K_{j-1}-(j-1)<k\leqslant K_{j-1}} \left\lfloor \frac{x}{k}-b\right\rfloor \\
+cx\sum_{K_{j+1}-(j+1)<k\leqslant K_{j+1}} F \left( \frac{x}{k}-b\right)-cx\sum_{K_{j-1}-(j-1)<k\leqslant K_{j-1}} F\left( \frac{x}{k}-b\right)\\
=\big((2j+1)\sqrt{j+1}-(2j-1)\sqrt{j-1}\big)\sqrt{cx}\\
+O\big((1+b)^2j^{5/2}(cx)^{-1/2}\big)+O\big((1+b)j\big).
\end{multline}
\end{prop}
\dem

Observons que, si $j=1$, la deuxième et la quatrième somme sont vides.

Pour $t\soe 1$, nous récrivons \eqref{180518b} sous la forme
\begin{equation}\label{180526b}
F(t)=\frac 1{\lfloor t\rfloor}-  \frac{a+b+1}{2t^2}+O\left((1+b)^2/t^3\right)
\end{equation}

D'après \eqref{180526c}, la contribution à la troisième et quatrième somme de \eqref{180526a} du terme d'erreur de~\eqref{180526b} est
\begin{equation}\label{180526d}
\ll j(1+b)^2(cx/j)^{-3/2}=j^{5/2}(1+b)^2(cx)^{-3/2}.
\end{equation}

On a ensuite, en utilisant \eqref{t45} et \eqref{180526c},
\begin{align}
\sum_{K_{j+1}-(j+1)<k\leqslant K_{j+1}} \frac{1}{\left(x/k-b\right)^2} &=\sum_{K_{j+1}-(j+1)<k\leqslant K_{j+1}} \Big((j+1)/cx+O\big((1+b)(j/cx)^{3/2}\big)\Big)\notag\\
&=(j+1)^2/cx+O\big(j^{5/2}(1+b)(cx)^{-3/2}\big),\label{180526e}
\end{align}
et, de même, 
\begin{equation}\label{180526f}
\sum_{K_{j-1}-(j-1)<k\leqslant K_{j-1}}\frac{1}{\left(x/k-b\right)^2}=(j-1)^2/cx+O\big(j^{5/2}(1+b)(cx)^{-3/2}\big).
\end{equation}

La proposition \ref{180522d} et \eqref{t45} nous donnent
\begin{multline}\label{180526g}
\sum_{K_{j+1}-(j+1)<k\leqslant K_{j+1}} \big(j\left\lfloor x/k-b\right\rfloor+cx/\left\lfloor x/k-b\right\rfloor\big) =\\
=\sum_{K_{j+1}-(j+1)<k\leqslant K_{j+1}} \big(2\sqrt{(j+1)cx}+O\big((1+b)^2j^{3/2}(cx)^{-1/2}\big)-\sum_{K_{j+1}-(j+1)<k\leqslant K_{j+1}} \left\lfloor x/k-b\right\rfloor\\
=(2j+1)\sqrt{(j+1)cx}+O\big((1+b)^2j^{5/2}(cx)^{-1/2}\big)+O\big((1+b)j\big)
\end{multline}
et, de même, 
\begin{multline}\label{180526h}
\sum_{K_{j-1}-(j-1)<k\leqslant K_{j-1}} \big(j\left\lfloor x/k-b\right\rfloor+cx/\left\lfloor x/k-b\right\rfloor\big) =\\
=(2j-1)\sqrt{(j-1)cx}+O\big((1+b)^2j^{5/2}(cx)^{-1/2}\big)+O\big((1+b)j\big)
\end{multline}

En regroupant les résultats \eqref{180526d}, \eqref{180526e}, \eqref{180526f}, \eqref{180526g} et \eqref{180526h}, on obtient le résultat annoncé. \fin

\smallskip

Passons maintenant à la contribution à $W_j$ des sommes de \eqref{180515d} où figure $a$.

\begin{prop}\label{180527b}
Pour $0<j\ioe cx/9(1+b)^2 -1$, on a
\begin{multline}\label{181002a}
j\sum_{K_{j-1}<k\leqslant K_{j+1}}\left( \left\lfloor \frac{x}{k}-a\right\rfloor -  \left\lfloor \frac{x}{k} -b\right\rfloor   \right)
+cx\sum_{K_{j-1}<k\leqslant K_{j+1}}\left( F\left( \frac{x}{k}-a\right) -  F\left( \frac{x}{k} -b\right)   \right)\\
=\frac{(2j-1)\sqrt{j+1}-(2j+1)\sqrt{j-1}}3\sqrt{cx}+R_j+O\big((1+b)^2j/c\big),
\end{multline}
avec
\begin{equation*}
R_j=R_j(x;a,b) = \sum_{K_{j-1}<k\leqslant K_{j+1}} (k^2c/x-j)\big(\{x/k-a\}-\{x/k-b\}\big).
\end{equation*}
\end{prop}
\dem

D'après la proposition \ref{180526i} et \eqref{180528c}, la contribution à la seconde somme du terme d'erreur de~\eqref{180518b} est
\[
\ll (K_{j+1}-K_{j-1})(1+b)^2(cx/j)^{-3/2} \ll (1+b)^2j/c^2x.
\]

Pour $k \ioe K_{j+1}$, l'inégalité \eqref{180526j} nous permet d'écrire
\[
 \frac{\{x/k-b\}-(a+b+1)/2}{(x/k-b)^2}=\frac{k^2}{x^2}\big(\{x/k-b\}-(a+b+1)/2\big) +O\big((1+b)^2k^3/x^3\big).
 \]
La contribution à cette seconde somme du terme $\big(\{t\}-(a+b+1)/2\big)/t^2$ de \eqref{180518b} est donc
\[
x^{-2}\sum_{K_{j-1}<k\leqslant K_{j+1}} k^2\big(\{x/k-a\}-\{x/k-b\}\big) +O\big((1+b)^2j/c^2x \big).
\]

En utilisant \eqref{180522b}, \eqref{180528c} et \eqref{180526j}, on voit que la contribution du terme $1/t$ de \eqref{180518b} est
\begin{align*}
\sum_{K_{j-1}<k\leqslant K_{j+1}} \Big(\frac 1{x/k-a}-\frac 1{x/k-b}\Big)&=-\frac c{x^2}\sum_{K_{j-1}<k\leqslant K_{j+1}} \big(k^2+O(bk^3/x)\big)\\
&=-\frac{(j+1)\sqrt{j+1}-(j-1)\sqrt{j-1}}{3\sqrt{cx}}+O(bj/cx).
\end{align*}

Enfin, la première somme du premier membre de \eqref{181002a} vaut, d'après \eqref{180522a},
\begin{multline*}
c(K_{j+1}-K_{j-1})-\sum_{K_{j-1}<k\leqslant K_{j+1}} \big(\{x/k-a\}-\{x/k-b\}\big)=\\
(\sqrt{j+1}-\sqrt{j-1})\sqrt{cx}+O(c)-\sum_{K_{j-1}<k\leqslant K_{j+1}} \big(\{x/k-a\}-\{x/k-b\}\big).
\end{multline*}

On obtient la relation \eqref{181002a} en collectant ces estimations.\fin

\smallskip

Notons que la majoration {\og triviale\fg} de $R_j$ est donnée par la proposition \ref{180528d} :
\begin{equation*}
R_j \ll (x/cj)^{1/2}+j.
\end{equation*}

\smallskip

Estimons enfin la contribution à $W_j$ des termes de \eqref{180515d} où figure $N_j$.

\begin{prop}\label{180527a}
Pour $0<j\ioe cx/9(1+b)^2 $, on a
\[
jN_j+cxF(N_j)=2\sqrt{jcx} +O\big((1+b)j\big).
\]
\end{prop}
\dem

Rappelons que $N_j$ est le plus grand nombre entier $n$ tel que
$$
(n+a)(n+b)\ioe cx/j.
$$
On a donc
\[
N_j=\left\lfloor \sqrt{cx/j+c^2/4}-(a+b)/2\right\rfloor.
\]

En particulier,
\[
\sqrt{cx/j}-(a+b)/2-1 \ioe N_j \ioe \sqrt{cx/j}.
\]

Pour $0<j\ioe cx/9(1+b)^2 $, la borne inférieure de cet encadrement de $N_j$ est $\soe \frac 23\sqrt{cx/j}$. Par suite, en utilisant \eqref{180519a},
\begin{align*}
jN_j+cxF(N_j)&=j\big(\sqrt{cx/j}+O(1+b)\big)+cx\Big(\frac{1}{\sqrt{cx/j}+O(1+b)}+O\big((1+b)/(cx/j)\big)\Big)\\
&=2\sqrt{jcx} +O\big((1+b)j\big).\fine
\end{align*}

\smallskip

Nous sommes maintenant en mesure d'adapter la proposition 6, p. 18 de \cite{zbMATH06754310}.

\begin{prop}\label{180528h}
Pour $0<j\ioe cx/9(1+b)^2-1$, on a
\[
W_j=f(j)\sqrt{cx} +R_j+O\big((1+b)^2j/c\big)+O\big((1+b)^2j^{5/2}(cx)^{-1/2}\big).
\]
avec
\begin{equation*}
f(j) =\frac{8j+2}3\sqrt{j+1}-\frac{8j-2}3\sqrt{j-1}-4\sqrt{j},
\end{equation*}
et où $R_j$ est défini par \eqref{181007b}.
\end{prop}
\dem

On obtient cette estimation en insérant dans \eqref{180515d} les résultats des propositions \ref{180527c}, \ref{180527b} et \ref{180527a}.~\fin
\section{Estimation de la somme $W(x;a,b)$}\label{180608c}

\subsection{Mise en évidence du terme principal}\label{180921a}

Nous donnons d'abord une expression de $W=W(x;a,b)$ faisant intervenir un paramètre réel~$J$.

\begin{prop}\label{180608a}
Pour $x$ et $J$ réels tels que $3/2 \ioe J \ioe cx/9(1+b)^2-1$, on a
\[
W=\frac{2}{\pi}\zeta(3/2)\sqrt{cx}+\Rcal(J)+O\big((1+b)^2J^2/c\big)+O\big((cx/J)^{1/2}\big)+O\big((1+b)^2J^{7/2}(cx)^{-1/2}\big),
\]
où l'on a posé
\[
\Rcal(J)=R_0+\sum_{1 \ioe j \ioe J} R_j,
\]
les quantités $R_j$ étant définies par \eqref{181007a} ($j=0$) et \eqref{181007b} ($j>0$).
\end{prop}
\dem


On a 
\begin{align*}
W &= W_0 +\sum_{1 \ioe j \ioe J}W_j + \sum_{j > J}W_j\\
&= W_0 +\sum_{1 \ioe j \ioe J}W_j +O\big((cx/J)^{1/2}\big)\expli{d'après la proposition \ref{180527d}, et car $1 \ioe cx/J$}\\
&=\frac 2{3}\sqrt{cx}+R_0+O\big((1+b)^2/c\big)\\
& \qquad \qquad+\sum_{1 \ioe j \ioe J}\Big(f(j)\sqrt{cx} +R_j+O\big((1+b)^2j/c\big)+O\big((1+b)^2j^{5/2}(cx)^{-1/2}\big) \Big)\\
&\qquad\qquad\qquad\qquad+O\big((cx/J)^{1/2}\big)\expli{d'après les propositions \ref{180528j} et \ref{180528h}}\\
&= \frac{2}{\pi}\zeta(3/2)\sqrt{cx}+\Rcal(J)+O\big((1+b)^2J^2/c\big)+O\big((cx/J)^{1/2}\big)+O\big((1+b)^2J^{7/2}(cx)^{-1/2}\big),
\end{align*}
où l'on a utilisé la somme de la série 
\[
\frac 2{3}+\sum_{j \soe 1}f(j)=\frac{2}{\pi}\zeta(3/2) \quad \expli{cf. \cite{zbMATH06754310}, (35), p. 20},
\]
et l'estimation $f(j) \ll j^{-3/2}$.\fin

\smallskip

En restreignant l'intervalle de variation de $J$, simplifions légèrement l'énoncé de la proposition~\ref{180608a}.
\begin{prop}\label{181017b}
Pour $x\soe 40(1+b)^3c^{-2}$ et $J$ tel que $3/2 \ioe J \ioe (x/c)^{1/3}$, on a
\[
W=\frac{2}{\pi}\zeta(3/2)\sqrt{cx}+\Rcal(J)+O\big((1+b)^2J^2/c\big)+O\big((cx/J)^{1/2}\big).
\]
\end{prop}
\dem

D'une part, on vérifie que
\begin{equation*}
x\soe 40(1+b)^3c^{-2} \implique (x/c)^{1/3} \ioe cx/9(1+b)^2 -1.
\end{equation*}

D'autre part, on a
\begin{equation*}
J \ioe (x/c)^{1/3} \ssi (1+b)^2J^{7/2}(cx)^{-1/2} \ioe (1+b)^2J^2/c.\fine
\end{equation*}

\subsection{Démonstration du Théorème A}\label{180608b}

\begin{prop}
Pour $x \soe 40(1+b)^4/c^3$,
\[
W=\frac{2}{\pi}\zeta(3/2)\sqrt{cx}+\Rcal(J)+O\big((1+b)^{2/5}c^{1/5}x^{2/5}\big),
\]
où $J= c^{3/5}(1+b)^{-4/5}x^{1/5}$.
\end{prop}
\dem

Nous choisissons $J$ pour équilibrer les deux premiers termes d'erreur de la proposition \ref{180608a} :
\[
J=c^{3/5}(1+b)^{-4/5}x^{1/5}.
\]

On vérifie que 
\[
x \soe 40(1+b)^4/c^3 \implique 3/2 \ioe J \ioe (x/c)^{1/3},
\]
et les deux termes d'erreur de la proposition \ref{181017b} sont
\begin{equation*}
\ll (1+b)^{2/5}c^{1/5}x^{2/5}.\fine
\end{equation*}

\section{Démonstration du Théorème B}\label{181020a}

\subsection{Rappels sur la théorie de van der Corput}

Afin d'estimer la somme $\Rcal(J)$, nous utiliserons l'énoncé suivant, dû à van der Corput, qui résulte de la version la plus simple de sa méthode d'estimation de sommes trigonométriques (cf.~\cite{zbMATH02600828}, Satz 5, p. 252 ; \cite{zbMATH02600829}, Satz 1, p. 215 ; \cite{zbMATH02599896}, (12), p. 23).
\begin{quote}
Soit $u,v$ des nombres réels tels que $v-u \soe 1$, et $f: [u,v] \vers \Real$ une fonction deux fois dérivable, dont la dérivée seconde est monotone et de signe constant. On a alors
\begin{equation}\label{180627a}
\sum_{u<n\ioe v} \left(\{f(n)\}-1/2\right) \ll \int_u^v\left\lvert f''(t)\right\rvert^{1/3} dt + \frac{1}{\sqrt{ \lvert f''(u)\rvert}}+\frac{1}{\sqrt{ \lvert f''(v)\rvert}}\virg
\end{equation}
où la constante implicite est absolue.
\end{quote}

Par sommation partielle, on déduit de \eqref{180627a} l'estimation suivante.
\begin{quote}
Soit $u,v$ des nombres entiers tels que $u\ioe v$, et $f: [u,v] \vers \Real$ une fonction deux fois dérivable, dont la dérivée seconde est monotone et de signe constant. Soit $(\gamma_n)_{u<n\ioe v}$ une suite de nombres réels. On a alors
\begin{equation}\label{180627c}
\sum_{u<n\ioe v} \gamma_n\left(\{f(n)\}-1/2\right) \ll G\int_u^v\left\lvert f''(t)\right\rvert^{1/3} dt + \frac{G}{\sqrt{ \lvert f''(u)\rvert}}+\frac{G}{\sqrt{ \lvert f''(v)\rvert}}\virg
\end{equation}
où 
\[
G=\max _{u<n\ioe v} \lvert \gamma_n\rvert +\sum_{u<n<v}\lvert \gamma_{n+1}-\gamma_n\rvert,
\]
et où la constante implicite est absolue.
\end{quote}

\subsection{Estimation de la somme $\Rcal(J)$}
\begin{prop}\label{180715b}
Pour $2 \ioe j \ioe (x/c)^{1/3}$, on a
\[
R_j \ll x^{1/3}j^{-1}+x^{1/4}j^{3/4}c^{-3/4}.
\]
\end{prop}
\dem

Pour $K_{j-1} <k \ioe K_{j+1}$, posons
\[
\gamma_k=\gamma_k(x,c,j)= ck^2/x-j.
\]

La suite $(\gamma_k)$ étant croissante, la proposition \ref{180711a} implique l'estimation
\[
\max_{K_{j-1} <k \ioe K_{j+1}} \lvert \gamma_n\rvert +\sum_{K_{j-1} <k < K_{j+1}}\lvert \gamma_{n+1}-\gamma_n\rvert \ll 1+j^{3/2}(c/x)^{1/2}\ll 1,
\]
en tenant compte de l'hypothèse $j \ioe (x/c)^{1/3}$.

En appliquant \eqref{180627c} aux deux fonctions $f(t)=x/t-a$ et $f(t)=x/t-b$, on obtient 
\begin{align*}
 \sum_{K_{j-1} <k \ioe K_{j+1}} \gamma_k\left ( \{x/k-a\}-\{x/k-b\}\right) &\ll 
\int_{K_{j-1}}^{K_{j+1}}(x/t^3)^{1/3} \, dt + \frac{1}{\sqrt{x/K_{j-1}^3}}+ \frac{1}{\sqrt{x/K_{j+1}^3}}\\
&\ll x^{1/3}\ln \frac{K_{j+1}}{K_{j-1}}+x^{-1/2}K_{j+1}^{3/2}\\
&\ll x^{1/3} \frac{K_{j+1}-K_{j-1}}{K_{j-1}}+x^{-1/2}(jx/c)^{3/4}\\
&\qquad\expli{d'après \eqref{180711b}}\\
&\ll x^{1/3}j^{-1}+x^{1/4}j^{3/4}c^{-3/4},
\end{align*}
où l'on a utilisé \eqref{180522a} et \eqref{180522z}. \fin

\begin{prop}\label{180715c}
Pour $x \soe 1+c$, on a
\[
R_0+ R_1 \ll x^{1/3}\ln (2x/c) +c^{-3/4}x^{1/4}.
\]
\end{prop}
\dem

On a
\[
R_0+R_1=\frac cx\sum_{1\ioe k\ioe K} k^2\big(\{x/k-a\}-\{x/k-b\}\big)+\sum_{0<k\leqslant K_{2}} (k^2c/x-1)\big(\{x/k-a\}-\{x/k-b\}\big).
\]

Un raisonnement analogue à celui de la démonstration de la proposition \ref{180715b} fournit l'estimation
\begin{align*}
R_0+R_1 & \ll 1+ \int_{1}^{K_{2}}(x/t^3)^{1/3} \, dt + \frac{1}{\sqrt{x}}+ \frac{1}{\sqrt{x/K_{2}^3}}\\
&\ll x^{1/3}\ln (2x/c) +c^{-3/4}x^{1/4}.\fine
\end{align*}

\begin{prop}\label{180715d}
Pour $x$ et $J$ réels tels que $3/2 \ioe J \ioe \big(x/(1+c)\big)^{1/3}$, on a
\[
\Rcal(J)\ll x^{1/3}\ln (x/c)+ x^{1/4}J^{7/4}c^{-3/4}.
\]

\end{prop}
\dem

D'après les propositions \ref{180715b} et \ref{180715c}, on a
\begin{align*}
\Rcal(J)&=\sum_{0 \ioe j \ioe J} R_j \\
&\ll x^{1/3}\ln (x/c) +\sum_{0 <j \ioe J}\big(x^{1/3}j^{-1}+x^{1/4}j^{3/4}c^{-3/4} \big)\\
&\ll x^{1/3}\ln (x/c)+ x^{1/4}J^{7/4}c^{-3/4}.\fine
\end{align*}

\subsection{Conclusion}\label{180627b}

\begin{prop}
Pour $x \soe 40 \, c^{-5}(1+b)^{27/2}$,
\[
W=\frac{2}{\pi}\zeta(3/2)\sqrt{cx}+O(c^{2/9}x^{4/9}).
\]
\end{prop}
\dem

Pour 
\begin{equation}\label{181018a}
x\soe 40(1+b)^3c^{-2}
\end{equation}
et $3/2 \ioe J \ioe \big(x/(1+c)\big)^{1/3}$, la conjonction des résultats des propositions \ref{181017b} et \ref{180715d} montre que 
\begin{equation*}
W-\frac{2}{\pi}\zeta(3/2)\sqrt{cx} \ll (cx/J)^{1/2} +x^{1/4}J^{7/4}c^{-3/4} +
x^{1/3}\ln (x/c)+(1+b)^2J^2/c.
\end{equation*}

Nous choisissons $J$ pour équilibrer les deux premiers termes d'erreur : $J=c^{5/9}x^{1/9}$. Cette quantité vérifie l'encadrement $3/2 \ioe J \ioe \big(x/(1+c)\big)^{1/3}$ si 
\begin{equation}\label{181018b}
x \soe 40 \max\big(c^{-5},(1+c)^{4}\big).
\end{equation}

On a alors, d'une part,
\[
(cx/J)^{1/2} +x^{1/4}J^{7/4}c^{-3/4}  \asymp c^{2/9}x^{4/9}.
\]
D'autre part, sous les hypothèses \eqref{181018a} et \eqref{181018b}, on a
\begin{equation*}
x^{1/3}\ln (x/c)\ll x^{1/3}(x/c)^{1/18} \ll c^{2/9}x^{4/9}, 
\end{equation*}
et l'inégalité 
\begin{equation}\label{181018c}
x\soe c^{-1/2}(1+b)^{9}
\end{equation}
entraîne
\begin{equation*}
(1+b)^2J^2/c\asymp (1+b)^2c^{1/9}x^{2/9} \ll c^{2/9}x^{4/9}, 
\end{equation*}

Si $x \soe 40 \,c^{-5}(1+b)^{27/2}$, les trois conditions \eqref{181018a}, \eqref{181018b} et \eqref{181018c} sont vérifiées, et le résultat est démontré.\fin

\medskip
\begin{center}
{\sc Remerciements}
\end{center}

{\footnotesize Le premier auteur remercie Julien Cassaigne de lui avoir suggéré de généraliser à $W(x;a,b)$ l'étude, effectuée dans \cite{zbMATH06754310}, du cas $a=1$, $b=2$.}

\providecommand{\bysame}{\leavevmode ---\ }
\providecommand{\og}{``}
\providecommand{\fg}{''}
\providecommand{\smfandname}{et}
\providecommand{\smfedsname}{\'eds.}
\providecommand{\smfedname}{\'ed.}
\providecommand{\smfmastersthesisname}{M\'emoire}
\providecommand{\smfphdthesisname}{Th\`ese}


\medskip

\footnotesize

\noindent BALAZARD, Michel\\
Institut de Math\'ematiques de Marseille\\
CNRS, Universit\'e d'Aix-Marseille\\
Campus de Luminy, Case 907\\
13288 Marseille Cedex 9\\
FRANCE\\
Adresse \'electronique : \texttt{balazard@math.cnrs.fr}

\bigskip

\noindent BENFERHAT, Leila\\
Institut de Math\'{e}matiques-USTHB\\
LA3C, Universit\'{e} des sciences et de la technologie Houari-Boum\'{e}di\`{e}ne\\
Bab Ezzouar\\
ALG\'{E}RIE\\
Adresse \'{e}lectronique: \texttt{lbenferhat@hotmail.com}

\bigskip

\noindent BOUDERBALA, Mihoub\\
Institut de Math\'{e}matiques-USTHB\\
LA3C, Universit\'{e} des sciences et de la technologie Houari-Boum\'{e}di\`{e}ne\\
Bab Ezzouar\\
ALG\'{E}RIE\\
Adresse \'{e}lectronique: \texttt{mihoub75bouder@gmail.com}

\end{document}